\documentclass{amsart}

\usepackage{amsthm,amssymb,verbatim}
\usepackage{graphicx}
\usepackage{enumerate}

\newcommand{\interior}{\operatorname{interior}}

 \newcommand{\MM}{\mathcal{M}}

 \newcommand{\RR}{\mathbf{R}}  
 \newcommand{\ZZ}{\mathbf{Z}}  
 \newcommand{\Sb}{\mathbf{S}}  
 \renewcommand{\SS}{\mathbf{S}}  

 \newcommand{\eps}{\epsilon}

 \newcommand{\genus}{\operatorname{genus}}

\input epsf
\def\begfig {
\begin{figure}
\small }
\def\endfig {
\normalsize
\end{figure}
}

\swapnumbers

    \newtheorem{theorem}    {Theorem}       [section]
    
    \newtheorem{corollary}  [theorem]     {Corollary}
    \newtheorem{proposition}       [theorem]       {Proposition}
    
    \newtheorem*{claim}{Claim}
    \newtheorem*{theorem*}{Theorem}

    \theoremstyle{definition}
    \newtheorem{definition}  [theorem] {Definition}

    \theoremstyle{definition}
    \newtheorem{remark}   [theorem]       {Remark}

\begin{document}

\title[]{On the number of minimal surfaces \\  with a given boundary}
\author{David Hoffman}

\address{Department of Mathematics\\ Stanford University\\ Stanford, CA 94305}
\email{hoffman@math.stanford.edu\\white@math.stanford.edu}
\author{Brian White}
\thanks{The research of the second author was supported by the NSF
  under grants~DMS-0406209 and DMS~0707126}

\date{July 6, 2008.}

\dedicatory{Dedicated to Jean  Pierre Bourguignon  on the occasion of his 60th birthday}
\begin{abstract}
We prove results allowing us to count, mod $2$, the number of embedded minimal surfaces of a specified topological type bounded by a curve $\Gamma \subset\partial N$, where $N$ is a mean convex $3$-manifold with piecewise smooth boundary. These results are extended to curves and minimal surfaces with prescribed  symmetries. The parity theorems are used in an essential manner to prove the existence of embedded genus-$g$ helicoids in $\SS ^2\times\RR$, and we give an outline of this application.
\end{abstract}

\maketitle

\renewcommand{\thesubsection}{\thetheorem}

\section{Introduction}

In \cite{tt1}, Tomi and Tromba used degree theory to solve a
longstanding problem about the existence of minimal surfaces with a
prescribed boundary: they proved that every smooth, embedded curve on
the boundary of a convex subset of $\RR ^3$ must bound an embedded
minimal disk.  Indeed, they proved that a generic such curve must
bound an odd number of minimal embedded disks. White \cite{wh4}
generalized their result by proving the following parity theorem.
 Suppose $N$ is a compact, strictly convex domain in $\RR^3$ with smooth boundary.
 Let $\Sigma$
be a compact $2$-manifold with boundary.   Then a generic smooth curve
$\Gamma\cong \partial\Sigma$ in $\partial N$ bounds an odd or even number
of embedded minimal surfaces  diffeomorphic to $\Sigma$ according to whether $\Sigma$
is or is not a union of disks.

In this paper, we generalize the parity theorem in several ways.
First, we prove (Theorem~\ref{DegreeTheoryStrict}) that the parity theorem holds for any compact riemannian $3$-manifold $N$
such that $N$ is strictly mean convex, $N$ is homeomorphic to a ball,  $\partial N$ is smooth, and
$N$ contains no closed minimal surfaces.
We then further relax
  the hypotheses by allowing $N$ to be  mean convex rather than 
strictly mean convex, and to have piecewise
  smooth boundary. Note that if $N$ is mean convex but not strictly 
mean convex, then $\Gamma$ might bound
  minimal surfaces that lie in $\partial N $.
We prove (Theorem~\ref{DegreeTheoryWeak}) that the parity theorem remains true for such $N$ provided
(1) unstable surfaces lying in $\partial N$ are not
counted, and (2) no two contiguous regions of $(\partial N)\setminus \Gamma$ are both smooth minimal surfaces.   We give examples showing that the theorem is false
without these provisos.

We extend the parity
theorem yet further  (see Theorem~\ref{DegreeTheoryWeakInvariant}) by showing that,
under an additional hypothesis,  it remains true for minimal surfaces
with prescribed symmetries.

The parity theorems described above are all mod $2$ versions of stronger results that describe
integer invariants.   The stronger results are given in section~\ref{IntegerInvariantSection}.

 The parity theorems are used in
an essential way to prove the
 the existence of
 embedded genus-$g$ helicoids in $\Sb^2\times \RR$.
 In Sections \ref{CountingHandles} and~\ref{Helicoid} we give a very brief outline of
 this application. (The full argument will appear in  \cite{hoffwhite3}.)

\section{Counting minimal surfaces}
\label{Counting}

Throughout the paper, $N$ will be a compact riemannian  $3$-manifold
and $\Sigma$ will be a fixed compact $2$ manifold.  If $\Gamma$ is an embedded curve
in $N$ diffeomorphic to $\partial \Sigma$, we let $\MM(N,\Gamma)$ denote the set of embedded minimal surfaces in $N$ that
are diffeomorphic to $\Sigma$ and that have boundary $\Gamma$.   We let
$|\MM(N, \Gamma)|$ denote the number of surfaces in $\MM(N,\Gamma)$.

In case $N$ has smooth boundary, we say that $N$ is  strictly mean convex 
provided the mean curvature is a (strictly) positive
multiple of the inward unit normal on a dense subset of $\partial N$.
\begin{theorem}\label{DegreeTheoryStrict}
Let $N$ be a smooth, compact, strictly mean convex
 riemannian $3$-manifold that is homeomorphic to
a ball and that has smooth boundary.
Suppose also
that $N$ contains no closed minimal surfaces.   Let $\Gamma\subset \partial N$
be a smooth curve diffeomorphic to $\partial \Sigma$.
 Assume that
 $\Gamma$ is bumpy in the sense that no  surface in
$\MM(N,\Gamma)$ supports a nontrivial normal Jacobi field with zero boundary values.

 Then $|\MM(N,\Gamma)|$ is even unless $\Sigma$ is a union of disks, in which case
 $|\MM(N,\Gamma)|$ is odd.
\end{theorem}

We remark that generic smooth curves $\Gamma\subset \partial N$ are
bumpy \cite{White-Indiana87}.

\begin{proof}
 Theorems 2.1 and 2.3 of \cite{wh4} are special cases of the theorem.
The proofs given there establish the more general result
here
provided one makes the following observations:
\begin{enumerate}
\item
There $N$ was assumed to be strictly convex, but exactly the same
proof works assuming strict mean convexity.
\item There $\Sigma$ was
assumed to be connected, but the same proof works for disconnected
$\Sigma$.
\item
In the proofs of Theorems~2.1  and~2.3 of \cite{wh4},
 the assumption
that $N$ is a subset of $\RR^3$  was used in order to invoke an
isoperimetric inequality, i.e., an inequality bounding the area of a
minimal surface in $N$ in terms of the length of its boundary. There
are compact mean convex $3$-manifolds for which no such isoperimetric
inequality holds.  However, if (as we are assuming here) $N$ contains
no closed minimal surfaces, then $N$ does admit such an isoperimetric
inequality \cite{wh6}.
 \item In the proofs in \cite{wh4}, one needs
to isotope any specified component of $\Gamma$ to a curve $C$ that
bounds exactly one minimal surface, namely an embedded disk. This was
achieved by choosing $C$ to be a planar curve.   For a general
ambient manifold $N$, ``planar" makes no sense.  However, any
sufficiently small, nearly circular curve $C\subset \partial N$
bounds exactly one embedded minimal disk and no other minimal
surfaces. (This property of such a curve $C$ is proved in the last
paragraph of \S3 in \cite{wh4}.)

\end{enumerate}
\end{proof}

\stepcounter{theorem}
\subsection{Mean convex ambient manifolds $N$ with piecewise smooth boundary}

For the remainder of the paper, we allow $\partial N$ to be piecewise smooth.  For simplicity, let
us take this to mean that $\partial N$ is a union of smooth $2$-manifolds with boundary (``faces''
  of $N$), any two of which are either disjoint or meet along a common edge with interior
  angle everywhere strictly between $0$ and $2\pi$.  (More generally, one could allow the
  faces of $N$ to have corners.)
We say that such an $N$ is  mean convex provided (1) at each interior point of each
face of $N$, the mean curvature vector is a nonnegative multiple of the inward-pointing unit
normal,  and
(2) where two faces meet along an edge, the interior angle is everywhere at most $\pi$.

The following example shows what can go wrong in Theorem~ \ref{DegreeTheoryStrict} if $N$ 
is mean convex
    but not strictly mean convex.
\\
\\
\noindent {\bf Example 1.} Let $N$ be a region in $\RR^3$ whose
boundary consists of an unstable catenoid $C$ bounded by two circles, 
together with the two disks bounded by those circles.   Note that $N$
is  mean convex with piecewise smooth boundary. Let $\Gamma$ be
a pair of horizontal circles in $C$ that are bumpy (in the sense of
Theorem~\ref{DegreeTheoryStrict}). Theorem~\ref{DegreeTheoryStrict} suggests that
$\Gamma$ should bound an even number of embedded minimal annuli in
$N$.
First consider the case when $\Gamma$ consists of two circles in $C$ very close to the waist circle. Then $\Gamma$  bounds precisely two minimal annuli. One of them is the  component of $C$ bounded by $\Gamma$. Because the circles in $\Gamma$ are close, this annulus is strictly stable. The other annulus bounded by $\Gamma$ is a strictly unstable catenoid lying in the interior of $N$. In order to get an even number of examples, we must count the stable catenoid lying on $C$.
Now suppose
the two components of $\Gamma$ are the two components of $\partial
C$. Then again $\Gamma$ bounds exactly two minimal annuli: the
unstable catenoid $C$, which is part of $\partial N$, and a strictly
stable catenoid that lies outside $N$.   Here, of course, we do not
count the stable catenoid since it does not lie in $N$.   Thus to get
an even number, we also must not count the unstable catenoid  that
lies in $\partial N$.

 This example motivates the following definition:

\begin{definition}
 $\MM^*(N, \Gamma)$ is the set of embedded minimal surfaces $M\subset N$ such that
 \begin{enumerate}
\item [i.)] $\partial M=\Gamma$,
\item [ii.)] $M$ is diffeomorphic to $\Sigma$, and
\item [iii.)] any connected component of $M$ lying in $\partial N$ must be stable.
\end{enumerate}
\end{definition}

Example 1
suggests that in order to generalize
Theorem~\ref{DegreeTheoryStrict}
to  mean convex $N$ with piecewise smooth boundary,
we should replace $\MM(N, \Gamma)$ by $\MM^*(N,\Gamma)$.   However,
even if one makes that replacement, the following
example shows that an additional hypothesis is required.
\\
\\
\noindent {\bf Example 2.} Let $N$ be a compact, convex region in
$\RR^3$ such that $\partial N$ is smooth and contains a planar disk
$D$.   Let $\Gamma$ be a pair of concentric circles lying in $D$.
Then $\Gamma$ bounds exactly one minimal annulus: the region in $D$
between the two components of $\Gamma$.  That annulus is strictly
stable and lies in $\partial N$. Thus $\Gamma$ is bumpy (in the sense
of Theorem~\ref{DegreeTheoryStrict}) and $|\MM^*(N,\Gamma)|=1$.
Consequently, if we wish $|\MM^*(N,\Gamma)|$ to be even (as
 Theorem~\ref{DegreeTheoryStrict}  suggests it should be), then we need an additional
hypothesis on $N$ and $\Gamma$.

 Note that in example 2, $(\partial N)\setminus \Gamma$ contains
 two contiguous connected components (a
planar annulus and a planar disk) both of which are minimal surfaces.
The additional hypothesis we require is that $(\partial N)\setminus
\Gamma$ contains no two such components.

\begin{theorem}\label{DegreeTheoryWeak}
Let $N$ be a smooth, compact,  mean convex  riemannian $3$-manifold that is homeomorphic to
a ball, that has piecewise smooth boundary, and that contains no closed minimal surfaces.
Let $\Gamma\subset \partial N$ be a smooth, embedded bumpy curve
diffeomorphic to $\partial \Sigma$.
Suppose that no two contiguous connected components of $(\partial N)\setminus \Gamma$
are both smooth minimal surfaces.

 Then $|\MM^*(N,\Gamma)|$ is
even unless $\Sigma$ is a union of disks, in which case $|\MM^*(N,\Gamma)|$  is odd.
\end{theorem}

\begin{proof}
Since $N$ is compact,   mean convex, and contains no closed
minimal surfaces, the areas of minimal surfaces in $N$ are bounded in
terms of the lengths of their boundaries \cite{wh6}.

If  $\partial N$ is smooth and has nowhere-vanishing mean curvature, the result
follows immediately from Theorem~\ref{DegreeTheoryStrict}. We reduce
the general case to this special case as follows. Note that we can
find a one-parameter family $N_t$, $0\le t\le \eps$, of mean convex
subregions of $N$ such that
\begin{enumerate}
\item [i.)]$N_0=N$, \item [ii.)]the boundaries $\partial N_t$ foliate
an open set containing $\partial N$. \item [iii.)] for $t>0$ small,
$\partial N_t$ is smooth and  the mean curvature of $\partial N_t$ is
nowhere zero and points into $N_t$.
\end{enumerate}
 For
example, we can let $\partial N_t$ be the result of letting $\partial
N$ flow for time $t$ by the mean curvature flow.

\begin{claim} Suppose $M_i$ are smooth embedded minimal surfaces in $N$ diffeomorphic to
   $\Sigma$
and that $\partial M_i\to \Gamma$ smoothly.  Then a subsequence of the $M_i$ converges smoothly to a limit $M\in \MM^*(N,\Gamma)$.
\end{claim}

\begin{proof}[Proof of claim]  By Theorem~3 in \cite{wh8},
a subsequence converges smoothly away from a finite set $S$ to a  limit surface $M$. The surface $M$ is smooth and embedded, though portions of it
may have multiplicity $>1$.  Indeed, the proof of Theorem~3 in \cite{wh8}
 shows that the multiplicity is $1$
and the convergence $M_i\to M$ is smooth everywhere unless an interior point of $M$ touches
  $\Gamma$.

In fact,  no interior point of $M$ can touch $\Gamma$. For suppose to
the contrary that the interior of $M$ touches $\Gamma$ at a point
$p$. Let $C$ be the connected component of $\Gamma$ containing $p$.
By the strong maximum principle, $M$ must contains a whole
neighborhood of $p\in \partial N$. Indeed, by the strong maximum
principle (or by unique continuation), $M$ must contain the two
connected components of $(\partial N)\setminus \Gamma$ on either side
of $C$. But by hypothesis, at most one of those components is a
minimal surface, a contradiction. This proves that no interior point
of $M$ touches $\Gamma$.

Consequently, as noted above, $M$ has multiplicity $1$ and the
convergence $M_i\to M$ is smooth everywhere.  Thus $M\in
\MM(N,\Gamma)$.

Now suppose some connected component $M'$ of $M$ lies in $\partial N$.  Then the
corresponding component $M_i'$ of $M_i$ converges smoothly to $M'$ from one side of $M$.
This one-sided convergence implies that $M'$ is stable.  Thus $M\in \MM^*(N, \Gamma)$.
This completes the proof of the claim.
\end{proof}

Continuing with the proof of Theorem~\ref{DegreeTheoryWeak}, note
that $\MM^*(N,\Gamma)$ is finite.  For if it contained an infinite
sequence of surfaces then, by the claim, it would contain a smoothly
convergent subsequence. The limit of that subsequence would be an
element of $\MM^*(N,\Gamma)$. But by bumpiness of $\Gamma$, the
elements of $\MM^*(N,\Gamma)$ are isolated.  The contradiction proves
that $\MM^*(N,\Gamma)$ is finite.

Let $\Gamma_t$, $0\le t\le \eps$, be a smooth one-parameter family of embedded curves
such that $\Gamma_0=\Gamma$ and such that $\Gamma_t\subset \partial N_t$.
\newcommand{\wN}{\widehat{N}}
Let $M_0^1, \dots M_0^k$ be the set of surfaces in $\MM^*(N, \Gamma)$.
By the implicit function theorem, we can (if $\eps$ is sufficiently small) extend these to
one-parameter families
\[
    M_t^i \in \MM^*(\wN, \Gamma_t)  \qquad (i=1,2,\dots, k; \,  0\le t\le \eps)
\]
where $\wN$ is a riemannian $3$-manifold containing $N$ in its
interior.

In fact, $M_t^i$ must lie in $N$ provided $\eps>0$ is chosen sufficiently small.
To see this, assume for simplicity that $\Sigma$ is connected.   If $M_0^i$ does not
lie in $\partial N$, then by the strong maximum principle, it is never tangent to $\partial N$, so
by continuity, $M_t^i \subset N$ for all sufficiently small $t$.
Now suppose that $M_0^i$ does lie in $\partial N$.
Then (by definition of $\MM^*(N, \Gamma)$) it is
strictly stable.   The strict stability implies that in fact $M_t^i$ lies in $N$ for sufficiently small $t$.

Indeed, $M_t^i$ must lie not only in $N$ but also in $N_t\subset N$,  for all sufficiently small $t$.  For let $T=T(t)\in [0,t]$ be the largest
number such that $M_t^i\subset N_T$.  If $T< t$, then $M_t^i$ would touch $\partial N_T$ at an
interior point, violating the maximum principle.  Hence $T=t$ and therefore $M_t^i\subset N_t$.

The claim implies that if $\eps$ is sufficiently small, then each surface in
$\MM^*(N_t,\Gamma_t)$ will be one of the surfaces in $M_t^1,\dots, M_t^k$.
We may also choose $\eps$  small  enough so that the $M_t^i$ all have nonzero nullity.
Then
\[
  |\MM^*(N, \Gamma)| = k = |\MM(N_t, \Gamma_t)|
\]
which must have the asserted parity by
Theorem~\ref{DegreeTheoryStrict} (applied to $N_t$ and $\Gamma_t$.)
\end{proof}

\stepcounter{theorem}
\subsection{Counting in the presence of symmetry}
In some situations, it is important to  be able to say something about the
number of minimal surfaces that are diffeomorphic to a specified surface $\Sigma$ and that possess specified symmetries.
Suppose  $G$ is a group of isometries of $N$.

\begin{definition} If $\Gamma$ is a $G$-invariant curve in $N$, we let $\MM_G^*(N, \Gamma)\subset  \MM^*(N,\Gamma)$ denote the set of surfaces
in $\MM^*(N,\Gamma)$ that are invariant under $G$.
 A boundary $\Gamma \subset \partial N$ is called $G$-bumpy if no surface in $\MM^*_G(N, \Gamma)$ has a nontrivial $G$-invariant
normal Jacobi field that vanishes on $\partial M$.
\end{definition}

Theorem~\ref{DegreeTheoryWeak}  has a natural extension to  $G$-invariant surfaces:

 \begin{theorem} \label{DegreeTheoryWeakInvariant}
Let $N$ be a smooth, compact,  mean convex  riemannian $3$-manifold that is homeomorphic to
a ball, that has piecewise smooth boundary, and that contains no closed minimal surfaces.
Let $G$ be a group of isometries of $N$.
Let $\Gamma\subset \partial N$ be a smooth curve that is $G$-invariant
and $G$-bumpy.
Suppose that no two contiguous components of $(\partial N) \setminus \Gamma$ are both minimal surfaces.
\newline
\noindent
Suppose also that
\begin{enumerate}
\item[$(*)$] $\Gamma=\partial \Omega$ for some $G$-invariant region $\Omega\subset \partial N$.
\end{enumerate}
Then $|\MM^*_G(N,\Gamma)|$ is
even unless $\Sigma$ is a union of disks, in which case $|\MM^*_G(N,\Gamma)|$  is odd.
\end{theorem}

\begin{remark}In Theorem~\ref{DegreeTheoryWeakInvariant}, the hypothesis that $N$ contains no closed minimal surfaces is equivalent to the hypothesis that $N$ contains no closed $G$-invariant minimal surfaces. See \cite{wh6},  Theorem~2.5.

\end{remark}

 \begin{proof}
 In general, the proof is exactly the same as the proof in the non-invariant case.
 However (see Observation (4) in the proof of Theorem~\ref{DegreeTheoryStrict}), to carry out the proof, one must be able to isotope the
 connected components of $\Gamma$ in a $G$-invariant way to arbitrarily small, nearly circular
 curves in $\partial N$.
 The hypothesis that $\Gamma=\partial \Omega$ for a $G$-invariant region
    $\Omega\subset\partial N$
ensures that such isotopy is possible.  (Indeed, it is equivalent to the existence of such $G$-invariant
isotopies.)
\end{proof}

We do not know whether Theorem~\ref{DegreeTheoryWeakInvariant} remains true without the hypothesis (*).

\section{An Integer Invariant}\label{IntegerInvariantSection}

Suppose $N\subset \RR^3$ is a compact, strictly convex set with smooth boundary.
In the introduction, we quoted Theorems~2.1 and~2.3 of \cite{wh4} as asserting that if
   $\Gamma\subset \partial N$ is a smooth,
bumpy curve diffeomorphic to $\partial \Sigma$, then
\begin{equation}\label{mod2}
    |\MM(N, \Gamma)|
    \cong
    \begin{cases}
    1  &\text{if $\Sigma$ is a union of disks, and} \\
    0 &\text{if not}
    \end{cases}
\end{equation}
where $\cong$ denotes congruence modulo $2$.

In fact, the conclusion in \cite{wh4} is actually much stronger than~\eqref{mod2}.
  To state that conclusion, we need
some terminology.

\newcommand{\MMeven}{\MM_\text{even}}
\newcommand{\MModd}{\MM_\text{odd}}

\begin{definition}
Let $\delta(\Sigma)=1$ if $\Sigma$ is a union of disks and $0$ if not.
If $\MM$ is a collection of smooth minimal surfaces, let
\[
    d(\MM) = |\MMeven| - |\MModd|
\]
where $\MMeven$ is the set of surfaces in $\MM$ with even index of instability and
$\MModd$ is the set of surfaces in $\MM$ with odd index of instability.
\end{definition}

With this terminology, the conclusion of Theorem~2.1 in \cite{wh4} is
\begin{equation}\label{integer}
    d(\MM(N, \Gamma)) = \delta(\Sigma).
\end{equation}
Note that~\eqref{integer} is stronger than~\eqref{mod2}.
  Indeed, \eqref{mod2} merely asserts that
   the two sides of~\eqref{integer} are congruent modulo $2$.
(See~\cite{Tromba1984} for a similar result for immersed minimal disks in $\RR^n$.)

If we start with the stronger conclusion~\eqref{integer}, then the arguments in \S2 produce stronger
versions of Theorems~\ref{DegreeTheoryStrict}, \ref{DegreeTheoryWeak}, ,
and \ref{DegreeTheoryWeakInvariant}:

\begin{theorem}
Under the hypotheses of Theorem ~\ref{DegreeTheoryStrict},
\[
    d(\MM(N,\Gamma))  = \delta(\Sigma).
\]
Under the hypotheses of Theorem  \ref{DegreeTheoryWeak},
\[
   d(\MM^*(N,\Gamma)) = \delta(\Sigma).
\]
Under the hypotheses of Theorem  \ref{DegreeTheoryWeakInvariant},
\[
    d_G(\MM^*_G(N, \Gamma)) = \delta(\Sigma)
\]
where $d_G(\cdot)$ is defined exactly like $d(\cdot)$, except that in
determining index of instability, we only count eigenfunctions that
are $G$-invariant.
\end{theorem}

The proofs are exactly as before.

%
%

\section{Counting the number of handles on a surface \\ invariant under an involution}
\label{CountingHandles}

Consider a minimal surface that has an axis of
 orientation preserving, $180^\circ$ rotational symmetry.
In many examples of interest, the handles of the surface are in some sense aligned along
the axis.   In this section, we make this notion precise, and we  observe that our parity
theorems apply to such surfaces.

Recall, for example, that Sherk constructed a singly periodic, properly embedded minimal surface
  $M\subset \RR^3$
 that is asymptotic to the planes $x=0$ and $z=0$ away from the $y$-axis, $Y$.
By scaling, we may assume that $M$ intersects $Y$ precisely at the
lattice points $(0,n,0)$, $n\in \ZZ$.
Now $M$ has various lines of orientation preserving, $180^\circ$ rotational symmetry.
For example, $Y$ is one such a line, and the line $L$ given by $x=z$, $y=1/2$ is another.
Intuitively, the handles of $M$ are lined up along  $Y$ but not along $L$.
(The surface $M$ is also invariant under $180^\circ$ rotation about the $x$ and $z$ axes,
but those rotations reverse orientation on $M$.)
We make the intuition into a precise notion by observing that the rotation about $Y$ acts
on the first homology group $H_1(M,\ZZ)$ by multiplication by $-1$, whereas rotation about $L$ acts
on $H_1(M,\ZZ)$ in a more complicated way.

\newcommand{\Surf}{S}   

\begin{proposition}\label{tTopology}
Suppose $\Surf$ is a noncompact $2$-dimensional riemannian manifold of finite topology.
Suppose that $\rho: \Surf \to \Surf$ is an
 orientation preserving
isometry of order two, and that $\Surf/\rho$ is connected.
Then the following are equivalent:
\begin{enumerate}
\item $\rho$ acts by multiplication by $-1$ on the first homology group $H_1(\Surf,\ZZ)$.
\item the quotient $\Surf/\rho$ is topologically a disk. \item
$\Surf$ has exactly $2-\chi(\Surf)$ fixed points of $\rho$,
where $\chi(\Surf)$ is the Euler
  characteristic of $\Surf$.
\end{enumerate}
\end{proposition}

\begin{corollary}\label{even-odd}
If the equivalent conditions (1)-(3) hold, then
the surface $\Surf$ has either one or two ends, according to
whether $\rho$ has an odd or even number of fixed points in
$\Surf$.
\end{corollary}

\begin{remark}\label{endscomment}
To apply Proposition~\ref{tTopology} and its corollary to a compact manifold $M$ with non-empty boundary, one
lets $\Surf=M\setminus \partial M$.  Of course the number of ends of $\Surf$ is equal to the
number of boundary components of $M$.
\end{remark}

\begin{proof}[Proof of Proposition~\ref{tTopology}]
Suppose that (1) holds.  Let $\pi: \Surf \to \Surf/\rho$ be the projection
and let $C$ be a closed curve in $\Surf/\rho$.  Then $C'=
\pi^{-1}(C)$ is a $\rho$-invariant cycle in $\Surf$ and thus (by (1)) it bounds a $2$-chain in $\Surf$.
Consequently $\pi(C') = 2C$ bounds a $2$-chain in $\Surf/\rho$.
Thus $2C$ is homologically trivial in $\Surf/\rho$.  But
$\Surf/\rho$ is orientable, so $H_1(\Surf,\ZZ)$ has no torsion.
Thus $C$ is homologically trivial in $\Surf /\rho$. Since $\Surf
/\rho$ is noncompact and connected with trivial first homology group, it
 must be a disk.  Hence (1) implies (2).

To see that (2) implies (1), suppose that (2) holds.
It suffices to show that any $\rho$-invariant $1$-cycle in $\Surf$ is
a boundary.  (For if $C_0$ is any cycle in $\Surf$, then
$C_0 + \rho(C_0)$ forms  a $\rho$-invariant cycle.) Since $\Surf$
is oriented, $H_1(\Surf,\ZZ)$ has no torsion, so it suffices to
show that any $\rho$-invariant cycle $1$-cycle in $\Surf$ must be a boundary mod~$2$.
Let $C\subset \Surf$ be any $\rho$-invariant closed curve, not
necessarily connected. We may assume that $C$ is smooth and in
general position, i.e., that the self-intersections are
transverse. By doing the obvious surgeries at the intersections,
we may assume in fact that $C$ is embedded.

Now $\pi(C)$ is a smooth, embedded, not necessarily connected,
closed curve in $\Surf/\rho$. Since $\Surf/\rho$ is
topologically a disk,  $\pi(C)$ bounds a region $\Omega$.
It follows that $C$ bounds the region $\pi^{-1}(\Omega)$.  Thus
$C$ is homologically trivial mod $2$. This completes the proof
that (2) implies (1).

Finally we show that (2) and (3) are equivalent.
Let $P$ be the number of fixed points of $\rho$. Consider a
triangulation of $\Surf/\rho$ such the fixed points of $\rho$ are
vertices in the triangulation, and consider the corresponding
triangulation of $\Surf$.   Then from Euler's formula one sees
that
\[
  \chi(\Surf) = 2\chi(\Surf /\rho) - P
  \]
or
\[
  P = 2\chi(\Surf/\rho) - \chi(\Surf).
\]
Thus $P=2-\chi(\Surf)$ if and only if $\chi(\Surf/\rho)=1$.  Since
$\Surf/\rho$ is orientable and connected, its Euler
characteristic is $1$ if and only if it is a disk.  This proves
that (2) and (3) are equivalent.
 \end{proof}

\begin{proof}[Proof of Corollary \ref{even-odd}]
Since $\Surf/\rho$ is a disk, it has exactly one end.  Since
$\Surf$ is a double cover of $\Surf/\rho$, it must have either one or
two ends. Since $\Surf$ is oriented,
\begin{equation}\label{chiS}
  \chi(\Surf) = 2c - 2g - e,
\end{equation}
where $c$ is the number of connected components, $g$ is the sum of
the genera of the connected components, and $e$ is the number of
ends.  Thus $e$ is congruent mod $2$ to $\chi(\Surf)$, which by
Proposition ~\ref{tTopology}  is congruent, mod $2$, to the number
of fixed points of $\rho$.
\end{proof}

\stepcounter{theorem}
\subsection{Counting $Y$-surfaces}\label{CountingYSurfaces}
Let $N$ be a riemannian $3$-manifold.
We suppose that $N$ has a geodesic $Y$ and an orientation preserving,
order two isometry $\rho=\rho_Y: N \to N$ for which the set of fixed points is $Y$.

\begin{definition}  Suppose $M \subset N$ is an orientable,   $\rho $-invariant surface
such that $\rho:M\to M$ preserves orientation and such that
 $(M\setminus \partial M)/\rho$ is connected.
We will say that   {\em $M$ is a $Y$-surface} if $\Surf := M\setminus\partial M$
 satisfies the equivalent conditions in Proposition~\ref{tTopology}.
\end{definition}

Suppose for example that $N=\RR^3$ and that $Y$ is a line.  Then $\rho=\rho_Y$ is $180^\circ$
rotation about $Y$.   If $M$ is a $\rho_Y$-invariant catenoid, then either $Y$ is the axis of rotational
symmetry of $M$, or else $Y$ intersects $M$ orthogonally at two points on the waist of $M$.
 In the first case, $\rho$ acts trivially on the first homology of $M$, so
$M$ is {\bf not} a $Y$-surface.  In the second case, $\rho$ acts by multiplication by $-1$
on the first homology of $M$, so $M$ is a $Y$-surface.

\begin{definition}  We let
\[
    \MM_Y^*(N, \Gamma) = \{ M\in \MM^*(N, \Gamma): \text{$M$ is a $Y$-surface} \}.
\]
We say that a curve $\Gamma\subset \partial N$ is {\em $Y$-bumpy}
 if no surface in $\MM_Y^*(N, \Gamma)$
 carries a nontrivial, $\rho_Y$-invariant, normal
   Jacobi field that vanishes on $\Gamma$.
\end{definition}

The following result is a version of Theorem~\ref{DegreeTheoryWeakInvariant}:

 \begin{theorem} \label{DegreeTheoryY} Let $N$ be a smooth, compact,  mean convex
 riemannian $3$-manifold that is homeomorphic to a ball, that has piecewise smooth boundary,
 and that contains no closed minimal surfaces.
 Suppose that $Y$ is a geodesic in $N$ and that $\rho=\rho_Y: N\to N$ is an orientation preserving,
 order two isometry of $N$ with fixed point set $Y$.

 Let $\Gamma\subset \partial N$ be a smooth, embedded, $\rho$-invariant, $Y$-bumpy
 curve that carries  a  $\rho$-invariant orientation.

 Suppose that no two contiguous components of $(\partial N) \setminus \Gamma$ are both minimal surfaces.

 Then $|\MM_Y^*(N,\Gamma)|$ is even unless $\Sigma$ is a union of disks, in which case
 $|\MM_Y^*(N,\Gamma)|$ is odd.
 \end{theorem}

 \begin{proof}
 The proof is almost identical to the proof of Theorem~\ref{DegreeTheoryWeakInvariant}.
  One lets the group $G$ in Theorem~\ref{DegreeTheoryWeakInvariant}
 be the group generated by $\rho$.  The hypothesis (*) there follows from the hypothesis here
 that $\Gamma$ carries a $\rho_Y$-invariant orientation.
 \end{proof}

 %
 %

\section{Higher genus helicoids in $\Sb ^2 \times \RR$}\label{Helicoid}
\stepcounter{theorem}
\subsection{A boundary value problem for minimal $Y$-surfaces}Our motivation in formulating Proposition~\ref{tTopology} and
Theorem~\ref{DegreeTheoryY} comes from the desire to construct
embedded minimal surfaces in $\Sb ^2 \times \RR$, each of  whose ends
is asymptotic to  a helicoid in $\Sb ^2 \times \RR$. Take as a model
of $\Sb ^2 \times \RR$ the space $\RR ^2\times \RR$ on which each
$\RR ^2\times \{t\}$ has the metric of the  sphere  pulled back by
inverse stereographic projection. (The radius of that sphere is fixed
but arbitrary.)  This model is missing a line, $Z^*=\{\infty\}\times
\RR$, which we append in a natural way to $\RR ^2\times \RR$ with the
aforementioned product metric. It is easy to verify that a standard
helicoid  $H\subset \RR ^3$ with axis $Z=\{(0,0,t)\,:\, t\in \RR\}$,
 an embedded and ruled surface,  is also a minimal surface in  $\Sb ^2 \times \RR$. Here, it has two axes, $Z$ and $Z^*$.
 By a slight abuse of notation, we will use $H$ to refer to this minimal surface in $\Sb ^2 \times \RR$ .

 The horizontal lines on the euclidean helicoid  are  great circles in  the totally geodesic level-spheres of $\Sb ^2\times \RR$,
 each circle passing through the antipodal points $({\bf 0} ,t)\in Z$ and  $({\bf\infty},t)\in Z^*$. Let  $$X=(\Sb ^2\times \{0\})\cap H,$$ and denote by $Y$
   the great circle  at height $0$ passing through  $O=({\bf 0},0)$, $O^*=({\bf \infty}, 0)$, and orthogonal to the great circle  $X$.
 Just as on the Euclidean helicoid, $\rho _{Y}$, order-two rotation about $Y$, is an orientation preserving
  involution of $H$. 
  Note that under our identification of $\Sb ^2\times \RR$
  with $\RR ^3$, each of the great circles on $H$ corresponds to a horizontal line passing throught the $z$-axis, and  the great circles $X$ and $Y$ are identified with the $x$- and $y-$axes of $\RR ^3$.

  Denote by $H^+$  one of the two components in the complement of $H$.
Then  for any $c>0$, $ \rho _Y$  is an orientation preserving
involution of the domain
 \begin{equation}\label{Nc}
 N_c= H^+\cap \{|z|<c\}.\end{equation}
 Note that $\partial N_c$ is  mean convex, consisting of three minimal surfaces:  $H\cap \{|z|<c\}$,
 and two totally geodesic hemispheres, $H^+\cap \{z=\pm c\}$ . We will label these minimal surfaces  $H_c$ and $S_{\pm c}$,  respectively.

 The set $H_c \setminus (Z\cup Z^*\cup X)$ has four components. Let $Q$ be the component whose boundary contains the three geodesics $X^+=\{(x,0,0)| x\geq 0\}$,\, $Z\cap \{ 0\leq z\leq c\}$, and $Z^*\cap \{ 0\leq z\leq c\}$. The ``quadrant'' $Q$ has a fourth boundary curve, which is one of the two semicircular components of  $\partial S_c \setminus (Z\cup Z^*)$. We label this semicircle $T_{c}$. Note that  $T_{-c}:=\rho_Y(T_c)$ lies in $\partial( \rho _Y(Q))$. 

  Fix a value of $c$ and let $N=N_c$. 
Consider the union $Q\cup \rho _Y (Q) $, and define $\Gamma \subset \partial N$ to be the boundary of
$Q\cup\rho _Y (Q)$. Then
\begin{equation}\label{Gamma}
\Gamma = (Z\cap H_c)\cup T_c\cup (Z^*\cap H_c)\cup T_{-c}\cup X.
\end{equation} 
See Figure~\ref{GammaFigure}.
 The first four segments of $\Gamma$ form a piecewise smooth curve with four corners. Adding the
 great circle $X$ produces a curve that is singular at $O=({\bf 0},0)$ and at $O^*=({\bf\infty},0)$, where there
 are right-angle crossings. Note that $\Gamma$ is $\rho _Y $-invariant.

 \begfig
  \vspace{.2in}
  \centerline{
\includegraphics[width=2.55in]{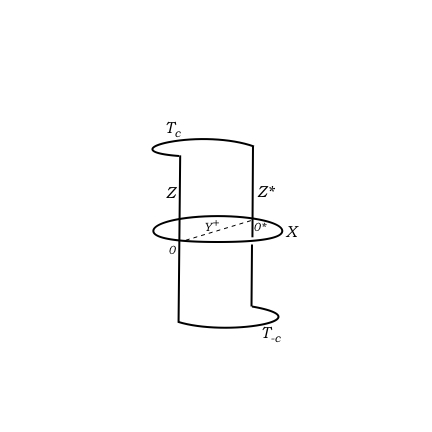}
   }
\begin{center}
 \vspace{-0.3in}
\caption{\label{GammaFigure}  The curve $\Gamma$.  In the figure, we have taken
$\RR^3 =\RR^2\times \RR$  as our model for $\Sb^2\times\RR$,  with the metric on 
 $\RR^2$  given by the pullback of the metric on $\Sb ^2$ via inverse stereographic projection. In this case, the pole of  $\Sb ^2$ is placed at the center of
  the semicircle $Y^-$. 
}.
 \end{center}
 \endfig
 If $\Gamma$ defined in \eqref{Gamma} is not $Y$-bumpy, we can make
 arbitrarily small perturbations of the curves $T_{\pm c}$ to
 make it so, while keeping the resulting curve in $\partial N$, and also $\rho_Y $-invariant.
  {\bf We will assume
 from now on that $\Gamma$ is $Y$-bumpy.}

 Suppose for the moment that we could produce a connected $Y$-surface $M \subset N$ with boundary
 $\Gamma$.  We will show in the next paragraph how this will enable us to construct
 a higher-genus helicoid.

 Since $\rho _Y |_M$ is orientation preserving, $Y$ must intersect $M$ orthogonally  in a discrete set of points, precisely
  the fixed points of  $\rho _Y |_M$. We will consider
 $M $ without its boundary, allowing us to apply Proposition~\ref{tTopology}. Namely, if $k=|Y\cap M|$, the number of points
 in $Y\cap M$, then

  $$k  =2- \chi (M ).$$
Extend $M$ by $\rho _Z$,  Schwarz reflection in $Z$ (or equivalently in $Z^*$), and let
\begin{equation}\label{Mtilde}
\tilde{M}=\interior (\overline{M \cup\rho _Z (M)}.
\end{equation}
  The surface $\tilde{M}$ is
smooth because $M$ is $\rho _Y$-symmetric,  and
$$|Y\cap \tilde{M}|=2k+2$$
because the points $O=({\bf 0},0)$ and at $O^*=({\bf \infty},0)$, which lie on $Y$, are in $\tilde{M}$. The surface $\tilde{M}$ is bounded
by  two great circles at levels $\pm{c}$. It is embedded because $\rho _Z (M)$ lies in $H^{-}$. Furthermore it
 is $\rho _Y  $-invariant by construction and satisfies the condition that $\rho_Y$ acts by multiplication by $-1$ on $H_1(M,Z)$.
 Therefore,
 $2k +2= 2-\chi(\tilde{M})$  by Proposition~\ref{tTopology}. Since $\tilde{M}$ has two ends, we have
$$2k +2 =2-(2-2\genus(\tilde{M}) -2),$$
or
$$\genus(\tilde{M})=k.$$
If we can produce $\tilde{M}=\tilde{M_c}$ for any cutoff height $c$, it is reasonable to expect that as $c\rightarrow\infty$, 
the $\tilde{M_c}$ converge subsequentially to  an embedded genus-$k$ minimal surface each of whose ends is asymptotic to $H$ or a rotation of $H$.
 In \cite{hoffwhite3}, we prove that this is the case.

\stepcounter{theorem}
\subsection{Existence of a suitable $M\in\MM _Y^*(N,\Gamma)$ with $|Y\cap M|=k$}
 How are we going to produce, for each positive integer $k$, a connected, embedded, minimal $Y$-surface $M \subset N$ with boundary
 $\Gamma$? The answer is: by induction on $k$, using Theorem~\ref{DegreeTheoryY}. The
 details, carried out in \cite{hoffwhite3}, are somewhat intricate. We describe
here the main idea and the intuition behind the proof.

  First of all, it would seem that Theorem~\ref{DegreeTheoryY} is not suited to prove
  existence of the desired surfaces because in most cases it asserts that the number of surfaces in a given class is even.
This could mean that there are {\em zero} surfaces in the class. We
begin to address this problem by dividing the class of surfaces
according to their geometric behavior near $O$. Why this helps
will be made clear below.

 Since we are working with one fixed domain, namely $N=N_c$ as defined in \eqref{Nc}, we will 
 suppress the reference to $N$ and write $\MM^*_Y(\Gamma)$ instead of $\MM^*_Y(N, \Gamma)$.
  We can decompose
$\MM_Y^*(\Gamma)$ into two sets by looking at how  a surface
$S\in\MM_Y^*(\Gamma)$ attaches to $\Gamma$ at the crossing $O$, the
intersection of  the vertical line $Z$ and the great circle $X$. The
line $Z$ is naturally oriented. Through $O$ also passes the great
circle $Y$, and our choice of the component $H^+$ of $\, \Sb ^2\times
\RR\setminus H$ lets us specify $Y^+:= Y\cap H^+$. We may now choose
a component of $X\setminus \{O,O^*\} $ to be $X^+$ by the condition
that at $O$,
 the tangent vectors to $X^+, Y^+$ and $Z^+$ form an oriented basis. The geodesics $X$,  $Z$, and $Z^*$ divide $H$ into four ``quadrants."
 We will call a quadrant whose boundary contains $Z^+\cup X^+$, or
$Z^-\cup X^-$ a {\em positive} quadrant, and refer to the other two quadrants as  {\em negative} quadrants.
\begin{definition} Given a nonnegative integer $k$,
\begin{itemize}
 \item[]$\MM_Y^*(\Gamma,k)\subset \MM^*_Y(\Gamma)$ is the collection
 of embedded minimal $Y$-surfaces  $M$ with the property that
 $|M\cap Y|=k$.
 \item[]$\MM_Y^*(\Gamma,k, +)\subset \MM_Y^*(\Gamma,k)$ is the subset  of surfaces tangent to the positive quadrants at $O$.
 \item[] $\MM_Y^*(\Gamma,k, -)\subset \MM_Y^*(\Gamma,k)$ is the subset of surfaces tangent to the negative
 quadrants at $O$.
 \end{itemize}
\end{definition}

Now we approximate $\Gamma$ by  smooth  embedded curves
$\Gamma (t)\subset\partial N$.  We have to do this in order to apply any of our parity theorems. We want the four corners to be rounded and the two crossings to be resolved. At $O$, we modify
$\Gamma$ in a small neighborhood of radius $t>0$ by connecting $Z^+$ to $X^+$ and $Z^-$ to $X^-$.
Given this choice at $O$, we resolve the crossing at $O^*$  according to whether $k$ is even or odd as follows: connect  positively if $k$ is even (i.e. $Z^+$ to $X^+$  and  $Z^-$ to $X^-$) and negatively (i.e. $Z^+$ to $X^-$  and  $Z^-$ to $X^+$)   if $k$ is odd.
Again we modify in a manner that preserves $\rho _Y$-invariance, and we choose $t$ small enough so that the neighborhoods of the corners and the crossings are pairwise-disjoint. We will refer to such a rounding as an {\em adapted positive
rounding of $\Gamma$}. Note that when $k$ is odd, an adapted positive rounding of $\Gamma$ is connected, while when $k$ is even,
such a rounding has two components. See Figure~\ref{roundings}.

\begfig
  \centerline{
\includegraphics[width=3.8in]{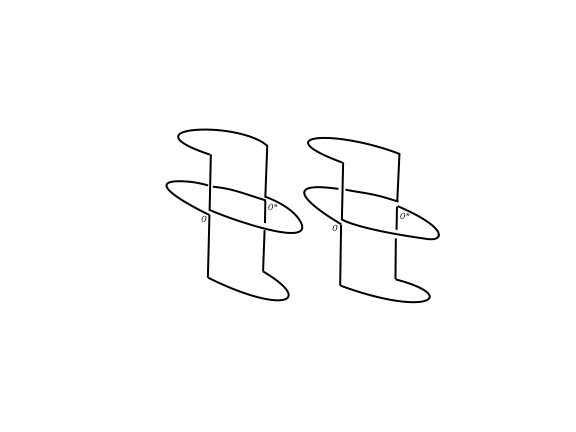}
   }
   \vspace{-0.3 in}
\begin{center}
 \vspace{-0.3in}
\caption{\label{roundings}  The two adapted positive roundings of $\Gamma$. On the left, the rounding at $O^*$  is the same as at the point $O$, resulting in a curve with two components. On the right,  the rounding at $O^*$ is positive to negative, resulting in a connected curve.
}.
 \end{center}
 \endfig
Our motivation for  the choice of desingularization at $O^*$ is given by the following
\begin{proposition} \label{kposneg}A surface $S\in \MM_Y^*(\Gamma,k, +)$ is tangent at
$O^*$ to the positive quadrants if $k$ is even, and to the negative
quadrants if $k$ is odd.
\end{proposition}
\begin{proof}
For any oriented surface $S$, we have  \eqref{chiS}  $$\chi(S)=2c(S)-2\,genus(S)-e(S),$$ where $e(s)$ is the number of ends of 
$S$,  $c(S)$ is  the number of components of $S$, and $genus(S)$ is the sum of the
genera of the components of $S$. If  $S\in\MM_Y^	*(\Gamma,k)$,  then using  Proposition~4.1 we have
\begin{equation}\label{k-e}
k=|Y\cap S| =2 -\chi(S)\cong e(S),
\end{equation}
where $\cong$ denotes equivalence mod $2$.
\begin{claim} If $S\in\MM(\Gamma,k,+)$, then
 $e(S)=\begin{cases}
2  & \text{if $S$ is positive at $O^*$,}\\ 
1  & \text{if $S$ is negative at $O^*$.}
\end{cases}
$
\end{claim} 
\noindent The proposition follows from the claim and the congruence  \eqref{k-e}.
\begin{proof}[Proof of Claim]Let $B(O)$ be a geodesic ball of radius $r>0$ centered at 
$O$, and let $B(O^{*})$ be the corresponding ball centered at $O^*$ with the same radius. We may choose $r$ small enough so that the surface $S^{'}=S \setminus (B(O)\cup B(O^{*}) )$ has the same number of ends as $S$: i.e., $e(S^{'})=e(S)$.  We may make 
$r$ smaller if necessary so that near $O$ (say in a geodesic ball of radius $2r$ centered at $O$),  the
boundary curve $\Gamma ^{'} =\partial S^{'}$ consists of a segment of $X^+$ joined to a segment of $Z^+$ by a single curve in $\partial B(O)$ together with a segment  
of $X^-$ joined to a segment of $Z^-$ by a single curve in $\partial B(O)$. It is precisely here that we have 
used the fact that $S\in \MM_Y^*(\Gamma,k, +)$ and not just in $\MM_Y^*(\Gamma,k)$.
Making $r$ smaller if necessary, we may assert that if $S$ is tangent to the positive quadrants at $O^{*}$, then
near $O^{*}$  the curve $\Gamma ^{'}$  connects positively, just  as  it does near $O$. This 
implies that $\Gamma ^{'}$ has two components. Therefore $e(S^{'})=2$. If  $S$ is tangent to the negative quadrants at $O^{*}$, then near $O^{*}$  the curve $\Gamma^{'}$ will connect $X^+$ to 
$Z^{-}$ and $X^-$ to $Z^+$. In this case,   $\Gamma ^{'}$ is  connected and $e(S^{'})=1$.
Since we chose $r$ small enough so that $e(S^{'})=e(S)$, we have proved the claim.
\end{proof}
\end{proof}

Let $\Gamma (t)$, $t>0$ small,  be a  smooth family of adapted positive roundings of $\Gamma$. We will round in such a way
that for each corner and crossing $q$,
$$\lim _{t\rightarrow 0} (1/t)(\Gamma (t)-q)$$
 is a smooth embedded curve, and such that $ \Gamma (t)$ converges smoothly to
 $\Gamma$ except perhaps at the corners and crossings of $\Gamma$. It is now reasonable to expect that if we specify a surface $M\in\MM_Y^*(\Gamma,k)$ as a sort
 of initial data at $\Gamma= \Gamma (0) $ we can deform it to a family of embedded
 minimal $Y$-surfaces $S_t\subset N$ with $\partial S_t =\Gamma (t)$. In fact we can do this in a unique manner.
 \begin{definition} For any nonegative integer $j$,   the set $\MM_Y^*(\Gamma (t), j)$  is the collection of embedded minimal $Y$-surfaces $S\subset N$ with $\partial S =\Gamma (t)$ and
 $|S\cap Y|=j$
 \end{definition}

  \begin{theorem}  Let $N= N_c\subset \Sb^2\times \RR$ be a domain of the form given in
 \eqref{Nc} for some fixed positive constant $c$. Let $\Gamma$ be the curve specified in
 \eqref{Gamma}, perturbed if necessary to become $Y$-bumpy.

   Let  $\Gamma (t)$, $t>0$ small,  be a  smooth family of adapted positive roundings of $\Gamma$.
 Suppose for some nonnegative integer $j$, that there
 exists a surface  $M \in\MM_Y^*(\Gamma,j)$. Then there exists a constant
 $a=a(\Gamma, M)>0$ such  that  for $ t<a$, each approximating curve $\Gamma (t)$ bounds an embedded minimal $Y$-surface $S_t$ with the following properties:

 \begin{enumerate}
 \item Each $S_t$ is the normal graph over a region $\Omega _t\subset\tilde{M}$ that is bounded by the projection of $\Gamma (t)$ onto $\tilde{M}$;
 \item  The family of surfaces $S_t$ is smooth in $t$ and  converges smoothly
 to $M $ as $t\rightarrow 0$;
 \item If $M \in\MM_Y^*(\Gamma , j,+)$, then $S_t\in \MM_Y^*(\Gamma (t) ,j)$, i.e.
 $|S_t\cap Y|=j$;
 \item If $M \in\MM_Y^*(\Gamma ,j,-)$, then $S_t\in \MM_Y^*(\Gamma (t) ,j+2)$,
 i.e.
 $|S_t\cap Y|=j+2$.
  \end{enumerate}

Furthermore, if $\hat{S}\in\MM_Y^*(\Gamma (t_0) ,j)$, $t_0<a$,  then it
lies in a  smooth one-parameter family of surfaces $S_t \in
\MM_Y^*(\Gamma (t),j)$, $t\leq t_0$, with the property that the family
has,  as a smooth limit as  $t\rightarrow 0$,  an embedded minimal
$Y$-surface $M \subset N$ that lies either in $\MM_Y^*(\Gamma, j)$
or in $\MM_Y^*(\Gamma , j-2)$.
 \end{theorem}

  Statements $(3)$ and $(4)$ have a simple geometric interpretation. Suppose we have a family of surfaces
 in $S_t\in\MM ^*_Y(\Gamma (t), k)$ for some smooth family $\Gamma (t)$ of adapted positive
 roundings of $\Gamma$. They will limit to an embedded minimal $Y$-surface $M\subset N$ with boundary $\Gamma$.  If they limit to an $M \in\MM_Y^*(\Gamma , j,+)$, then
 the points $S_t\cap Y$ stay bounded away from the crossings $\{O,O^*\}$.  Hence the
 $S_t$  have the property that $|S_t\cap Y|=|M\cap Y|=j$. However, if they limit to an
 $M \in\MM_Y^*(\Gamma , j,-)$, then each of the $S_t$ is a graph over a region $\Omega _t$ that contains both $O$ and $O^*$. Two points are lost. Hence  $j=|S_t\cap Y|=|M\cap Y|+2$.

 The theorem above tells us that there is a correspondence between surfaces
 in $\MM_Y^* (\Gamma (t), k)$ and   surfaces in the union of 
 $\MM^*_Y(\Gamma ,k,+)$ with $ \MM^*_Y(\Gamma ,k-2,-)$:
   \begin{corollary}\label{inductionstep}
 $$|\MM^*_Y(\Gamma (t),k)|=|\MM^*_Y(\Gamma ,k,+)|+|\MM^*_Y(\Gamma ,k-2,-)|.$$
 \end{corollary}

 We can now carry out the induction.
 Again,  $\cong$ denotes congruence modulo $2$. In our situation, the number of ends of a surface  $S\in\MM _Y^*(\Gamma, k)$ is one or two, so the number of components of $S$ is at most two. (See Corollary~\ref{even-odd} and Remark~\ref{endscomment}.) Since $S$ is a $Y$-surface  we know, by
 Proposition~\ref{tTopology},   that $k=|S\cap Y| =2-\chi(S)$.  It is easy to see that when
  $k=1$ (or $k=0$),  $S$  is a disk (or the union of two disks).
Corollary~\ref{inductionstep} and Theorem~\ref{DegreeTheoryY} yield in this  situation that
 $$1\cong |\MM^*_Y(\Gamma (t),k)|=|\MM^*_Y(\Gamma ,k,+)|+|\MM^*_Y(\Gamma, k-2,-)|=|\MM^*_Y(\Gamma ,k,+)|,$$
 the last equality being simply the fact that it is impossible for a surface to intersect $Y$ in a
 negative number of points.  Therefore we have established the existence of the desired surface  for
 $k=0$ or $k=1$.  In fact we get existence of a surface in $\MM^*_Y(\Gamma ,k,+)$.
 However there is nothing special in this context about being in
  $\MM^*_Y(\Gamma ,k,+)$ as opposed to being in $\MM^*_Y(\Gamma ,k,-)$. If we redid the entire construction  by starting out by requiring our smoothing to be negative at
  $O$, we would wind up with an odd number of surfaces in  $\MM^*_Y(\Gamma ,k,-)$, for $k=0$
  and $k=1$.

 Now assume $k\geq 2$, and suppose that for any $j<k$, that $|\MM^*_Y(\Gamma ,j,+)|\cong\MM^*_Y(\Gamma ,j, -)\cong 1$. Corollary~\ref{inductionstep} together with Theorem~\ref{DegreeTheoryY} yield in our situation that
$$0\cong |\MM^*_Y(\Gamma (t),k)|=|\MM^*_Y(\Gamma ,k,+)|+|\MM^*_Y(\Gamma , k-2,-)|.$$
 But $|\MM^*_Y(\Gamma ,k-2,-)|\cong 1$, by assumption. Therefore
  $0\cong|\MM^*_Y(\Gamma ,k,+)|+1,$
 or $$ |\MM^*_Y(\Gamma ,k,+)|\cong 1.$$  Hence,  this class of surfaces is not empty for any  nonnegative integer $k$.  As indicated above, the same is true for
 $\MM^*_Y(\Gamma ,k,-)$.
   Whether or not we have produced two geometrically {\em different} (i.e. non-congruent)
  solutions to our problem turns out to depend on whether
  $k$ is even or odd---but that is another story.
\nocite{hoffwhite1} \nocite{hoffwhite2}
\bibliography{genusone}
\bibliographystyle{alpha}

 \end{document}